\documentclass[12pt,notitlepage,twoside]{article}
\usepackage{showkeys}
\usepackage{latexsym} 
\usepackage{amsmath} 
\usepackage{amsfonts} 
\usepackage{amssymb}
\usepackage{graphicx}
\newcommand{\ristretto}{{\setbox0 =\hbox{$\mathsurround=0pt\vert$}\hbox
  {\lower\dp0 \copy0 }}}

\renewcommand{\thesection}{\arabic{section}.} 
\renewcommand{\le}{\leqslant}
\renewcommand{\ge}{\geqslant}
\newcommand{\abs}[1]{\left\vert #1 \right\vert}
\newcommand{\eps}{\varepsilon}
\newcommand{\imp}{\Rightarrow}
\newcommand{\W}{\mathcal W}
\newcommand{\hh}{\mathcal H}
\newcommand{\s}{\subseteq}
\newcommand{\fhi}{\varphi}

\newcommand{\U}{\mathcal U}
\newcommand{\R}{\mathbb R}
\newcommand{\N}{\mathbb N}
\newtheorem{prop}{Proposition}[section]
\newtheorem{lemma}{Lemma}[section]
\newtheorem{dfn}{Definition}[section]
\newtheorem{theorem}{Theorem}[section]
\newtheorem{corollary}{Corollary}[section]

\newcounter{rema}[section]
\newcounter{exa}[section]
\newenvironment{proof}{\noindent{\textbf{Proof.}}}{$\hfill\square$\vspace{0.1 cm}\\}
\newenvironment{example}{\stepcounter{exa}\noindent{\textbf{Example \thesection\arabic{exa}.}}}{$\hfill\blacktriangleleft$\\}
\newenvironment{rmk}{\stepcounter{rema}\noindent{\textbf{Remark \thesection\arabic{rema}.}}}{$\hfill\lhd$\\}
\usepackage{amssymb}
\begin{document}

\author{{\sc Mauro Garavello}
 \thanks{E-mail: \texttt{mgarav@sissa.it}. } }
\title{\textbf{Verification Theorems for Hamilton-Jacobi-Bellman equations}}

\date{SISSA-ISAS\\Via Beirut, 2-4\\34014 Trieste, Italy\\[15pt]
April 2002}
\maketitle
{\footnotesize 
\begin{abstract} 
 
\noindent 
We study an optimal control problem in Bolza form and we consider the 
value function associated to this problem. We prove two verification 
theorems which ensure that, if a function $W$ satisfies some suitable 
weak continuity assumptions and a 
Hamilton-Jacobi-Bellman inequality outside a countably $\mathcal H^n$-rectifiable set, then it is lower or equal to 
the value function. These results can be used for optimal synthesis
approach.



\end{abstract}

}
%
%
%
%
\section{Introduction.}
In this paper we consider a control system of the type:
\begin{equation}\label{eqintro}
\dot x=f(t,x,u), \qquad u\in U
\end{equation}
where $x\in \R^n$ is the state, $U\subset \R^q$ is the control space and
$f$ is the controlled dynamic. Given a target $S\subset \R^n$, a running
cost $L(t,x,u)$, a final cost $\psi(t,x)$ and an initial condition
$(t_0,x_0)$, we consider the optimal control problem in Bolza form consisting in
minimizing the integral of $L$ summed with the value of $\psi$ at
final points for trajectories that start at $x_0$ at time $t_0$ and reach the target
$S$. We define in the usual way the value function $V(t_0,x_0)$ to be the
infimum of the problem with initial condition $(t_0,x_0)$.  It is well known
that, under special conditions, $V$ satisfies the
Hamilton-Jacobi-Bellman equation in viscosity sense \cite{bardi} and
it is the unique solution.
Part of the proof is based on the Dynamic Programming Principle.

Therefore given a function $W$ with suitable properties, it is possible to determine if $W$
coincide with the value function, checking if it is a viscosity solution
to the HJB equation. This type of theorems, called verification theorems,
are useful, for example, when a candidate value function is produced
by means of the construction of a synthesis \cite{ps}.
It is then natural to ask for minimal conditions under which
a function $W$ coincides with the value function.
If we know that $W$ was obtained via a synthesis then the inequality
$W\geq V$ is granted by construction, thus we take this assumption.
Then, for $W$ to coincide with the value function, we prove it is sufficient that, outside a rectifiable
set of codimension one, both $W$ is differentiable and it
satisfies a Hamilton-Jacobi-Bellman inequality in classical sense.
Moreover, we make use of only some weak continuity assumptions, already used in \cite{ps} to prove optimality of a regular extremal synthesis, 
see Theorem \ref{teo1} and Theorem \ref{teo2} for details. A first result in this direction can be found in
\cite{F-R}, where the HJB inequality is asked outside
a locally finite collection of regular manifolds of positive codimension
(under more restrictive continuity assumptions). Notice that, for an optimal control problem, if the value function is also semiconcave, it is differentiable outside a countably $\mathcal H^n$-rectifiable set, see \cite{cms}.

We start considering the main assumptions for the problem and presenting two
technical lemmas, one of which dealing with the cardinality of the
intersections between admissible trajectories and a countably $\mathcal
H^n$-rectifiable set, while the other giving some conditions to assure the
monotonicity of a real valued function. Also we state, without proofs, two
propositions dealing with the properties of the solution to (\ref{eqintro})
and in particular dealing with existence, uniqueness and continuous dependence
by data.

Then, in Section~\ref{se6}\!\!, we recall briefly the synthesis approach
and various results available in the literature for comparison.
Some examples of regular optimal synthesis, to which our main results
are applicable, are given.

The first case we treat is the problem of finite time. We define a value
function as the infimum, over all admissible trajectories reaching the target
in finite time. The main result of this part is Theorem \ref{teo1} which
permits to verify if the function $W$ is lower or equal than the value
function.

Next, we consider the infinite time problem. In this case the value function
(\ref{vfit}) is defined as the infimum of the cost functional over all
admissible trajectories reaching the target in infinite time. The main result
of this section is Theorem \ref{teo2} which gives sufficient conditions on the
function $W$ to ensure the inequality $W\le V$, where $V$ is the value
function. In this case, for a technical reason, we consider a suitable
neighborhood $S_1$ of the target $S$ and we suppose that the final cost $\psi$
is defined on $S_1$ in order to give sense to the limit in the definition of
the value function (\ref{vfit}). As a corollary of Theorem \ref{teo1} and
Theorem \ref{teo2} we can treat a mixed case (see also \cite{p}), considering
at the same time the trajectories reaching the target both in finite time and
in infinite time.

A key ingredient for Theorem \ref{teo1} and Theorem \ref{teo2} is the
positiveness of the Lagrangian $L$, in order to prevent some bad phenomena
such as the permanence of the system for an arbitrary interval of times in a
region where $L$ is negative making the value function equal to $-\infty$ as
we see in Example \ref{ExaLpos}1. More precisely, it is not necessary to
suppose $L$ positive in the whole space, but some relaxed assumptions can be
taken, as we see in Remark \ref{RemL}4.

This paper ends with an appendix, where we give the definition of a non
continuous viscosity solution as in \cite{bardi} and we state Theorem
\ref{teoA}, which ensures that, under suitable assumptions, the value
functions (\ref{vf}) and (\ref{vfit}) are viscosity solutions to the
Hamilton-Jacobi-Bellman equation.

\begin{center}

{\bf Acknowledgments}

\end{center}
The author wishes to thank Prof. B. Piccoli, for having proposed him the study of this problem and for his useful advice, and the referees, for their improving suggestions.\vspace{20pt}

%
%
%
%
\setcounter{equation}{0}
\section{Preliminaries.}\label{se2}
We consider a control system:
\begin{equation}\label{cs}
\dot x(t)=f(t,x(t),u(t)),\quad (t,x)\in\Omega,\quad u(t)\in U
\end{equation}
where
\renewcommand{\labelenumi}{(A-\theenumi )}
\begin{enumerate}
\item $\Omega$ is an open and connected subset of $\R\times\R^n$.\label{a1}
\item $U$ is a non-empty subset of $\R^q$, for some $q\ge 1$, $q\in\N$.\label{a2}
\item $\U =L^p(\R ;U)$ with $1\le p<+\infty$ is the set of admissible controls.\label{a3}
\item \label{a4} $f:\Omega\times U\to\R^n$ is measurable in $t$, continuous in $(x,u)$, differentiable in $x$ and, for each $u\in U$, $D_xf(\cdot,\cdot,u)$ is bounded on compact sets. Moreover there exists $\fhi_1:\R\to\R^+$ integrable and for every $K$, compact subset 
of $\Omega$, there exist a modulus of continuity 
$\omega_K$ and a constant $L_K>0$ such that, if $(t,x)\in K$ and $(t,y)\in K$, then for all $u$
\begin{equation}\label{fc}
\begin{cases}
\abs{f(t,x,u)-f(t,y,u)}\le\omega_K(\abs{x-y})\\
\left(f(t,x,u)-f(t,y,u)\right)\cdot (x-y)\le L_K\abs{x-y}^2\\
\abs{f(t,x,u)}\le L_K(\fhi_1(t)+\abs{u}^p).
\end{cases}
\end{equation}
\end{enumerate}
We consider a function $L:\Omega\times U\to\R$ and assume:
\begin{enumerate}
\item[(A-5)] $L$ is measurable in $t$ and continuous
in $(x,u)$. Moreover, there exist
$\fhi_2 :\R\to\R^+$ integrable and, for every $R\ge 0$, $C_R\ge 0$ such 
that
\begin{equation}\label{lc}
\abs{L(t,x,u)}\le C_R(\fhi_2(t)+\abs{u}^p),\quad \abs{(t,x)}\le R.
\end{equation}
\end{enumerate}
In this paper we indicate with $x(\,\cdot\,;u,t_0,x_0)$ the solution to
(\ref{cs}) such that $x(t_0;u,t_0,x_0)=x_0$. Define the value function:
\begin{equation}\label{vf}
V(t_0,x_0):=\!\!\!\inf_{\substack{u\in\U\vspace{3pt}\\ (T,x(T;u,t_0,x_0))\in S} }\!\!\!\left\{\int_{t_0}^T\!\!\!\!\!\!L(s,x(s;u,t_0,x_0),u(s))ds\!+\!\psi (T,x(T;u,t_0,x_0))\right\}
\end{equation}
where $S$ - the target - is a closed subset of $\R\times\R^n$ contained in
$\Omega$, $\psi :S\to\R$ is the final cost. We recall the 
following
definition:\\\vspace{-.4cm}

\begin{dfn}
A subset $A$ of $\Omega$ is a \textsl{countably $\mathcal H^n$-rectifiable set} if there exist $A_1$ and $A_2$ such that $A=A_1\cup A_2$, $A_1$ is a finite or countable union of connected $\mathcal C^1$ submanifolds of positive codimension, and $\hh^n(A_2)=0$, where $\hh^k$ is the $k$-dimensional Hausdorff measure.
\end{dfn}
%
%
%
%
\section{Examples of syntheses.}\label{se6}
\setcounter{equation}{0}
In next sections we give sufficient conditions for a candidate value
function $W$ to coincide with $V$. Beside some regularity conditions, we 
ask a HJB inequality outside a countably $\mathcal H^n$-rectifiable set. 
This regularity is shared by every function $W$ obtained from
a regular synthesis, thus it can 
be used to prove the optimality of the synthesis itself. 
In this section we give various examples to which Theorem \ref{teo1} is 
applicable.
First of all, we need some definitions.
\begin{dfn}
A synthesis $\Gamma$ is a collection $\{(x_{(\bar t,\bar y)}(\cdot),
 u_{(\bar t,\bar y)})\}_{(\bar t, \bar y)\in \Omega}$ such that
$x_{(\bar t,\bar y)}(\cdot)=x(\cdot;u_{(\bar t,\bar  y)},\bar t,\bar 
y):[\bar t,\tau(\bar t,\bar y)]\to\R^n$,  
$u_{(\bar t,\bar y)}\in\U$ for every $(\bar t,\bar y)\in\Omega$,
$x_{(\bar t,\bar y)}(\tau(\bar t,\bar y))\in S$  and  
for every $t\in[\bar t,\tau(\bar t,\bar y)]$  
\[
u_{(t,x_{(\bar t,\bar y)}(t))}(s)=u_{(\bar 
t,\bar y)}(s+t) \quad\textrm{a.e.}
\]
and 
\[
x_{(t,x_{(\bar t,\bar y)}(t))}(\cdot) 
=x_{(\bar t,\bar y)}(\cdot+t)
\]
\end{dfn}
\begin{dfn}
A synthesis $\Gamma$ is optimal if every $u_{(\bar t,\bar y)}$ is an optimal
control. 
\end{dfn}

There is a standard method in geometric control theory to construct an 
optimal synthesis, see \cite{BP-98}. This consists of four steps: 1) 
using Pontryagin  Maximum Principle and other geometric tools to study the 
properties of optimal trajectories, 2) derive a sufficient family of 
extremal trajectories (i.e. trajectories satisfying PMP), 3) construct
a synthesis formed by extremal trajectories and 4) prove its optimality.
In many cases, for autonomous systems, it happens that the extremal 
synthesis is associated to a feedback $u:\R^n\to U$ that is smooth on
each stratum of a stratification, see \cite{ps} for details.
Roughly speaking a stratification is a locally finite collection of
disjoint regular submanifolds, of various dimensions, that is a partition
and such that the boundary of each manifold is union of manifolds of
higher codimensions. In this case the synthesis is called regular in the 
sense of Boltyanskii-Brunovsk\'y, see \cite{Bolt,Brun,ps}.

Step 4) of the geometric control approach can thus be obtained in 
essentially two ways: either using the regularity of the synthesis, see
\cite{ps}, or proving that the candidate value function $W$ associated to 
the synthesis coincides with $V$. The latter is exploited in \cite{F-R}
for a continuous $W$, defined on a subset of $\R^n$,
that is differentiable and satisfies the 
HJB equation outside a locally finite union
of smooth submanifolds of positive codimension. Then the optimality
is granted for initial points for which all admissible trajectories
remains in the domain of $W$.
A mild generalization is obtained in \cite{Bressan}, where trajectories 
can exit the domain of $W$, but the boundary of the domain of $W$ is a 
level set of $W$ itself. 
Another approach is the one of nonsmooth analysis, 
using which various verification theorems can be proved, see 
for example \cite{Vinter}.

Our main results, see Theorems~\ref{teo1} and \ref{teo2},
generalize previous results in the following way:

\renewcommand{\labelenumi}{\theenumi. }
\begin{enumerate}
\item As in \cite{Bressan} we assume that 
$W$ can be defined on a subset and the boundary of its domain is a level
curve of $W$.
\item We ask $W$ to be differentiable and satisfy HJB only outside
a countably ${\mathcal H}^n$-rectifiable set.
\item $W$ is only lower semicontinuous (satisfying other weak continuity
assumptions).
\end{enumerate}
A direct comparison with results of nonsmooth analysis is difficult.
However, we point out that the value function fails in general to be
locally Lipschitz continuous, see Example 3.1, for regular synthesis.
In case of locally Lipschitz regularity, our result is consequence of
those obtained by nonsmooth analysis methods, see 
for example \cite{Clarke,Vinter}.

We give now some examples to illustrate the applicability of our results.
A whole class of examples can be find in \cite{BressanPiccoli,bpic}.
The first example shows a typical regular synthesis with a non locally
Lipschitz continuous value function. In the second, the value function is not
continuous and it is differentiable only outside a countably ${\mathcal
H}^n$-rectifiable set. Last example shows the well known Fuller
phenomenon. In this case optimal trajectories have an infinite number
of switchings and the methods of Boltyanskii-Brunovsk\'y do not work
(while it does the result of \cite{ps}).\vspace{.5cm}

\begin{example}
Let $x\in\R$ and $u\in[-1,1]$. Consider the control system
\[
\ddot x+x=u
\]
and the problem of reaching the origin in minimum time. If 
we 
define $x_1=x$ and $x_2=\dot x$ we obtain the following first-order 
system:
\begin{equation}\label{sin1}
\left\{
\begin{array}{l}
\dot x_1=x_2\\
\dot x_2=-x_1+u.
\end{array}
\right.
\end{equation}
Every optimal trajectory is a bang-bang trajectory, i.e. formed by arcs 
corresponding to control $+1$ or $-1$.
The synthesis is illustrated in Figure \ref{fig2}. 
There are some "switching curves":
\begin{itemize}
\item all semi-circles of radius $1$ contained in $\{(x_1,x_2):x_2\le0\}$ and centered at $(2n+1,0)$, with $n\in\N\setminus\{0\}$;
\item all semi-circles of radius $1$ contained in $\{(x_1,x_2):x_2\ge0\}$ and centered at $(-2n-1,0)$, with $n\in\N\setminus\{0\}$.
\end{itemize}
Optimal trajectories switch along these curves, i.e. change control from 
$+1$ to $-1$ or viceversa.
Let $\gamma^{\pm}$ be the trajectory that switches at points $(\pm 2,0)$ 
(defined say on $[-\infty,0]$). Then the 
value function is not locally Lipschitz continuous at any point of supp 
$(\gamma^\pm)$,
but however it satisfies all the hypotheses of Theorem \ref{teo1}.
\begin{figure}
\begin{center}
\includegraphics[height=6cm]{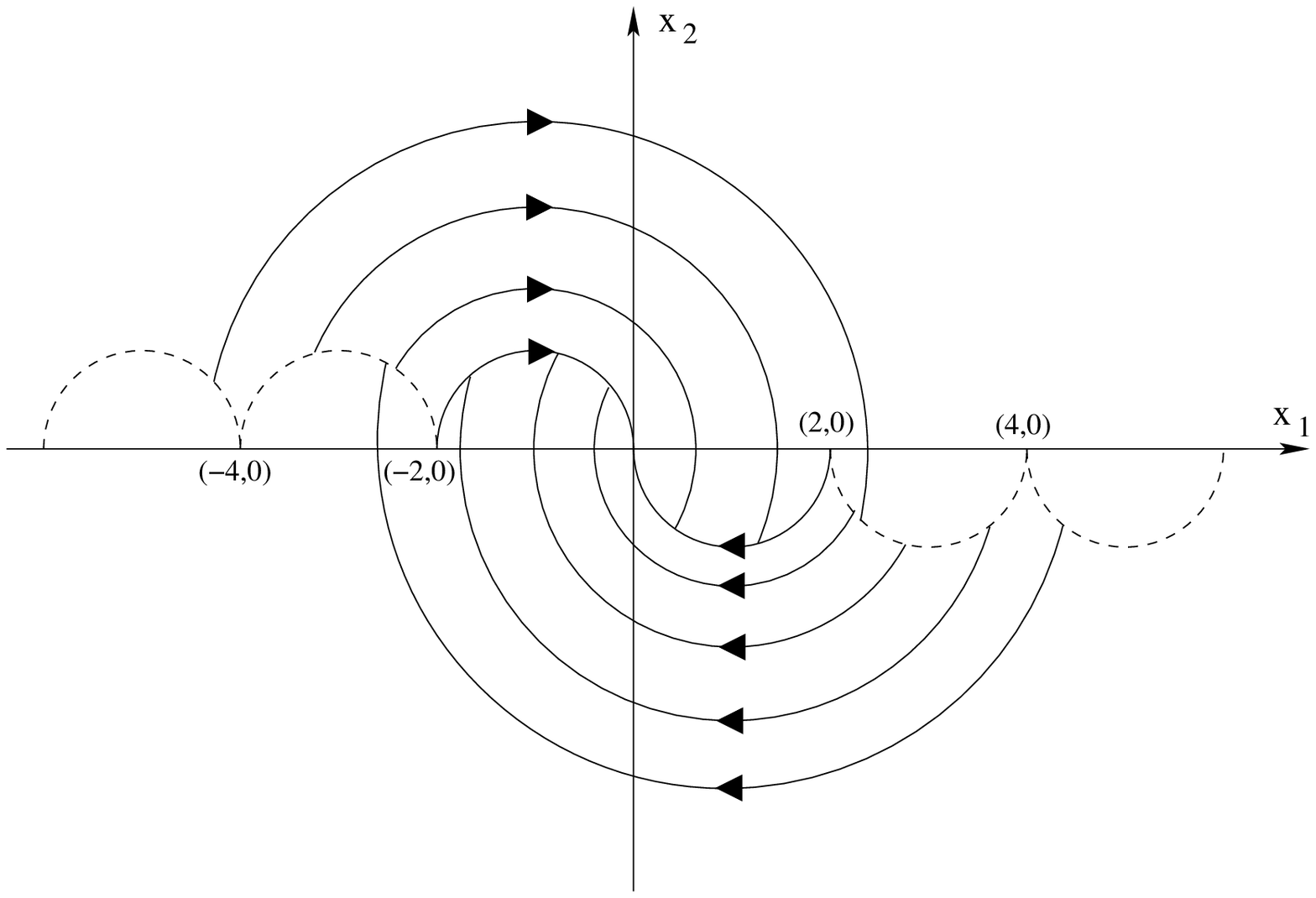}
\caption{Synthesis of system (\ref{sin1}).}
\label{fig2}
\end{center}
\end{figure}
\end{example}

\begin{example}
Let $\Omega=\R^2$, $f\equiv0$, $L\equiv1$. Consider the target:
\[
S=\left\{(t,x):x\ne0, t=\sin(1/x)\right\}\cup\left\{x=0,-1\le t\le1\right\}\cup\left\{t\ge1\right\}
\]
and the final cost $\psi$ constantly equal to $0$. The value function for this problem is given by:
\[
V(t,x)=\left\{
\begin{array}{lrl}
\sin(1/x)-t & \textrm{if} & x\ne0,t\le\sin(1/x)\\
1-t & \textrm{if} & x\ne0,\sin(1/x)<t<1\\
0 & \textrm{if} & t\ge1\\
-1-t & \textrm{if} &x=0,t\le-1\\
0 & \textrm{if} &x=0,-1<t<1.
\end{array}
\right.
\]
This function satisfies all the hypotheses of Theorem \ref{teo1} and clearly it is not continuous. Moreover it is differentiable outside a countably $\mathcal H^n$-rectifiable set $A$, which is not a locally finite union of regular manifolds.
\end{example}

\begin{example} 
\textbf{(Fuller phenomenon)}. Let us consider the system
\[ 
\left\{ \begin{array}{l} \dot x_1=x_2\\ \dot x_2=u \end{array} \right.
\] with
\begin{figure} 
\begin{center} 
\includegraphics[height=6cm]{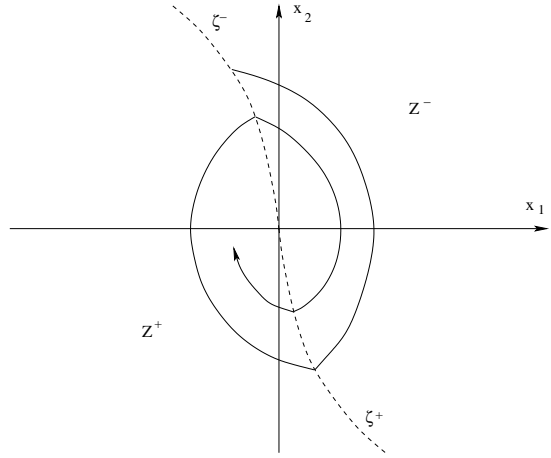}
\caption{Synthesis of Fuller phenomenon.} \label{Ful.Ph.} 
\end{center}
\end{figure} 
$\abs{u}\le1$, $\Omega=\R\times\R^2$, $S=\R\times\{0\}$,
$\psi\equiv0$, $L(t,x_1,x_2,u)=x_1^2$. This problem is well-known in the
literature, see for example \cite{Z-B}. Every optimal trajectory is
composed by an infinite number of bang-bang arcs, while the time for
reaching the origin of $\R^2$ is finite. There are two switching curves
$\zeta^+$ and $\zeta^-$ which separate $\R^2$ into two regions
$Z^+$ and $Z^-$ where the optimal trajectory uses respectively
the control $u=+1$ and $u=-1$, see Figure \ref{Ful.Ph.}. The value
function of this problem satisfies all the hypotheses of Theorem 
\ref{teo1}. 
\end{example} 

%
%
%
%
\section{Some useful results.}\label{se3}
\setcounter{equation}{0} 
We start by recalling without proofs some
classical results about ODEs. 
\begin{prop} \textbf{(Local existence and uniqueness of the trajectory).}
Assume (A-1)-(A-4). Fixed $u\in\U$ and $(t_0,x_0)\in\Omega$, there exist
$\delta>0$ and a unique absolutely continuous function
$x(\cdot;u,t_0,x_0):[t_0,t_0+\delta]\to\R^n$ solution to (\ref{cs}).
\end{prop}  
\begin{prop} \textbf{(Continuous dependence by
data).} Assume (A-1)-(A-4). Let $(t_0,x_0)\in\Omega$, $(t_0,x_n)\in\Omega$
for every $n\in\N$ and $u\in\U$, $u_n\in\U$ for every $n\in\N$. Let us
suppose that there exists a time $T>t_0$ such that $x(\cdot;u,t_0,x_0)$
and $x(\cdot;u_n,t_0,x_n)$ are defined in $[t_0,T]$. If $x_n\to x_0$ and
$u_n\to u$ in the strong topology of $L^p([t_0,T];U)$ as $n\to+\infty$,
then $x(\cdot;u_n,t_0,x_n)\to x(\cdot;u,t_0,x_0)$ uniformly in $[t_0,T]$
as $n\to+\infty$.  
\end{prop} 
Now, we present two technical lemmas 
used to prove the theorems of the next sections.  
%
%
%
%
\begin{lemma}\label{fin} Fix an element $\omega\in U$, $t'<t''$
and $x\in\R^n$ with $(t'',x)\in\Omega$. Assume that there exists $\W$, an
open neighborhood of $x$ in $\R^n$, such that $\zeta^y(\cdot)$, the
solution to $\dot\zeta^y(t)=f(t,\zeta^y(t),\omega )$ with
$\zeta^y(t'')=y$, is defined on $[t', t'']$ for any $y\in\W$ and
$(t,\zeta^y(t))\in\Omega\,\,\,\,\forall t\in[t',t'']$. Let $A$ be a
countable $\mathcal H^n$-rectifiable set.\\ Then for a.e. $y\in\W$ the set
$B^y:=\left\{t\in [t',t'']:(t,\zeta^y(t))\in A\right\}$ is finite or
countable.
\end{lemma}

This lemma is a slight generalization of a result proved in Theorem 2.14 of
\cite{ps}, since here we consider the trajectory coupled with
time.\vspace{.5cm}

\begin{proof}
We can write $A=A_1\cup A_2$, where $A_1=\cup_jM_j$ and $\{M_j\}_{j\in J}$ is a finite or countable family of connected submanifolds of $\R^{n+1}$ of codimension $d_j>0$, and $\hh^n (A_2)=0$. After replacing each $M_j$ by a finite or countable family of open submanifolds of $M_j$, we may assume that the $M_j$ are embedded. Define $\widetilde \W :=]t',t''[\times\W$ and let $\Phi$ be the map $\widetilde \W\ni (t,y)\mapsto (t,\zeta^y(t))\in \Omega$. The Jacobian of $\Phi$ is 
\begin{equation}\label{jphi}
{\bf J}\Phi\, =\,\left(\!
\begin{array}{l|c}
1&0\,\,\cdots\cdots\,\, 0\\
\hline {\bf b} & \begin{matrix}
                           \, & & \,\\[.1pt]
                           & {\bf V} ^\zeta (t;t',{\bf Id})&\\
                           & &
			   \end{matrix} 
\end{array}\right)\!
\end{equation}
where ${\bf b}$ is the column vector $f(t,\zeta^y(t),\omega)$ and ${\bf V}^\zeta (t;t', {\bf Id})$ is the fundamental matrix solution to the linear system
\begin{equation}\label{lin}
\dot v(t)=-D_xf(t,\zeta^y (-t+t'+t''),w)\cdot v(t)
\end{equation}
such that ${\bf V}^\zeta (t';t',{\bf Id})={\bf Id}$. So the determinant of
${\bf J}\Phi$ is equal to the determinant of ${\bf V}^\zeta (t;t',{\bf
Id})$, which is equal to
$exp\int_{t'}^ttr(-D_xf(s,\zeta^y(-s+t'+t''),\omega ))ds$, by Liouville's
theorem (see \cite{hart}). In particular $det ({\bf J}\Phi )$ is strictly
positive for any $t\in [t',t'']$. Moreover, by (A-\ref{a4}) $tr(-D_xf)$ is
bounded on compact sets and then there exist $c>0$, $C>0$ such that $0<c\le 
det ({\bf J}\Phi)\le C$.\\
So $\Phi$ is a Lipschitz diffeomorphism. In particular we have
$\hh^n(\Phi^{-1}(A_2))=0$. Now for each $j$ consider $\widetilde M_j
:=\Phi^{-1}(M_j)$. It is an embedded submanifold of codimension $d_j>0$.
Let $\Pi :\widetilde\W\to\W$ be the canonical projection. Consider the set
$S_j$ consisting of the points $s\in\widetilde M_j$ such that
$\Pi\ristretto_{\widetilde M_j}$ is not regular. Thus, by Sard's theorem,
$\mathcal L^n(\Pi (S_j))=0$. Moreover $\hh^n(\Pi(\Phi^{-1}(A_2)))=0$. So
the set $\mathcal B:=\Pi(\Phi^{-1}(A_2))\cup (\bigcup_j\Pi (S_j))$ has
Lebesgue measure $0$ in $\R^n$.\\
Let $y\in\W\setminus\mathcal B$. Then $(t,\zeta^y(t))\not\in A_2$ if $t'<t<t''$. To obtain the thesis, it is sufficient to show that, for each $j$, the set $E_j=\left\{ t\in ]t',t''[: (t,\zeta^y(t))\in M_j\right\}$ is at most countable. Fix $j$ and suppose $t\in E_j$. $\widetilde M_j$ has codimension $d_j>0$, so the dimension $\nu_j$ of $\widetilde M_j$ is less or equal to $n$. Since $y\not\in\mathcal B$, the map $d\Pi (t,y):T_{(t,y)}\widetilde M_j\to\R^n$ is onto, thus $\nu_j =n$ and $d\Pi (t,y)$ is injective. Obviously $d\Pi (t,y)(\frac{\partial}{\partial t})=0$, so $\frac{\partial}{\partial t}\not\in T_{(t,y)}\widetilde M_j$ and, consequently, $(\tilde t,y)\not\in\widetilde M_j$ if $0<\abs{\tilde t-t}\le\eps$ for $\eps >0$ sufficiently small. Therefore $t$ is an isolated point of $E_j$ and so the lemma is proved.
\end{proof}
%
%
%
%
\vspace{-.5cm}
\renewcommand{\theenumi}{\alph{enumi}}
\renewcommand{\labelenumi}{(\theenumi )}
\begin{lemma}\label{lerv}
Let $g$ be a real-valued function on a compact interval $[a,b]$. Assume that there exists a finite or countable subset $E$ of $[a,b]$ with the following properties:
\begin{enumerate}
\item $\liminf_{h\downarrow 0}\frac{g(x+h)-g(x)}{h}\ge0$ for all $x\in[a,b[\setminus E$,\label{gc1}
\item $\liminf_{h\downarrow 0}g(x+h)\ge g(x)$ for all $x\in[a,b[$,\label{gc2}
\item $\liminf_{h\downarrow 0}g(x-h)\le g(x)$ for all $x\in]a,b]$.\label{gc3}
\end{enumerate}
Then $g(b)\ge g(a)$.
\end{lemma}
For a proof of this lemma see \cite[Lemma B.1]{ps}.

%
%
%
%
\section{Problem with finite time.}\label{se4}
\setcounter{equation}{0}
We indicate with $\partial Q$ the topological boundary of an arbitrary $Q\s\R\times\R^n$.
Before stating the theorem we need the following definition

\begin{dfn}
Suppose that we have a time-varying Lipschitz-continuous vector field $X$
on $\R^n$ and
$W:\Omega\to\R\cup\{\pm\infty\}$. We say that $W$ has the \textsl{no
downward jumps} property (NDJ) along $X$ if
for any $[a,b]\ni t\mapsto\gamma(t)$, solution to
$\dot\gamma(t)=X(t,\gamma(t))$ such that $(t,\gamma(t))\in\Omega$ $\forall
t\in [a,b]$, we have $\liminf_{h\downarrow 0} W(t-h,\gamma(t-h))\le
W(t,\gamma(t))$, whenever $t\in ]a,b]$.
\end{dfn}  

\begin{theorem}\label{teo1}
Suppose (A-\ref{a1})-(A-5) hold. Let $Q\s\Omega$ be an open subset 
containing $S$. Let $W:\overline Q\to\R$ be a lower semicontinuous function
such that:
\renewcommand{\theenumi}{\roman{enumi}} 
\renewcommand{\labelenumi}{\theenumi)} 
\begin{enumerate} 
\item \label{NU1} $W$ has the NDJ property along every time-varying
vector field  of the type $f(t,x,u)$ with $u\in U$ fixed and for each $t$
\[
{\textrm{ess-liminf}}_{y\to x}W(t,y)\le W(t,x).
\]
\item $W\le\psi$ on $S$.\label{i2}
\item \label{i3} At every point $(t,x)\in\partial Q$ one has
\[
W(t,x)=\sup_{(s,y)\in Q}W(s,y).
\]
\item\label{i4} There exists a countably $\mathcal H^n$-rectifiable set
 $A\s\Omega$ such that $W$ is differentiable on $Q\setminus A$ and satisfies
\[
W_s(s,y)+\inf_{\omega\in U}\left\{ W_y(s,y)\cdot
f(s,y,\omega)+L(s,y,\omega)\right\}\ge 0\quad \textrm{on }Q\setminus A.
\]
\item $L\ge 0$.\label{i5} 
\end{enumerate}
Then $W\le V$ on $Q$. If $Q=\Omega$ we can drop hypotheses \ref{i3}) and \ref{i5}). 
\end{theorem}
\begin{proof}
Suppose by contradiction that there exists
$(t_0,x_0)\in Q$ such that $W(t_0,x_0)>V(t_0,x_0)$. In particular $V(t_0,x_0)<+\infty$. First of all, let us consider the case $V(t_0,x_0)>-\infty$. So we can find
$\eps>0$, $\delta>0$ such that
\begin{equation}\label{1}
V(t_0,x_0)\le W(t_0,x_0)-2\eps
\end{equation}
and, by the lower semicontinuity of $W$,
\begin{equation}\label{2}
\abs{x-x_0}<\delta \quad\imp\quad W(t_0,x)>V(t_0,x_0)+\eps.
\end{equation}
We can find $u^\ast\in\U$ such that $x^\ast (\cdot ):=x(\cdot; u^\ast,t_0,x_0)$ satisfies $(T,x^\ast (T))\in S$ and
\begin{equation}\label{3}
\int_{t_0}^TL(s,x^\ast (s), u^\ast (s))ds+\psi (T,x^\ast (T))\le V(t_0,x_0)+\frac{\eps}{2}\,\,.
\end{equation}
Moreover, for every $l\in\N$ there exists $u_l\in\U$ such that $\| u_l-u^\ast\|_{L^p([t_0,T])}\le\frac{1}{l}$, $u_l$ piecewise constant and left continuous. By \cite[Th\'eor\`em IV.9]{brez}, there exists a subsequence of $(u_l)_l$, denoted again by $(u_l)_l$, and a function $h\in L^p([t_0,T])$ such that $\vert u_l\vert\le h$ a.e. and $u_l$ converges to $u^\ast$ a.e. as $l\to+\infty$.
Hence, if we
denote by $x_l(\cdot)$ the trajectory $x(\cdot;u_l,T,x^\ast(T))$, 
for $l$ sufficiently big, we have (see Proposition 4.2),
\begin{equation}\label{5}
\abs{x_l (t)-x^\ast (t)}<\frac{\delta}{2}\quad\quad\forall t\in [t_0,T]
\end{equation}
and 
\begin{equation}\label{4}
\abs{\int_{t_0}^T[L(s,x_l(s),u_l(s))-L(s,x^\ast (s),u^\ast (s))]ds}\le\frac{\eps}{2}.
\end{equation}
Fix $l$ such that (\ref{5}) and (\ref{4}) hold and an interval $]t',t'']$ such
that $u_l (t)\equiv\omega$ on $]t',t'']$. Suppose that $(t,x_l (t))\in
Q\,\,\,\forall t\in [t',t'']$. Let $\zeta^y(t)$ be the trajectory associated
to the constant control $\omega$ such that $\zeta ^y(t'')=y$. By the fact
that $d(\partial Q,\{ (t,x_l (t)):t\in [t',t'']\} )>0$, we can find an open
neighborhood $\W$ of $x_l(t'')$ in $\R^n$ such that $(t'',y)\in 
Q\,\,\,\forall y\in\W$ and
$\left\{ (t,\zeta^y(t)):t\in[t',t'']\right\}\s Q\,\,\,\forall y\in\W$. By
Lemma \ref{fin}, we have that for a.e. $y\in\W$ the set 
$B^y:=\left\{t\in [t',t'']:(t,\zeta^y(t))\in A\right\}$
is at most countable.\\ 
Therefore,since for every fixed $t$ $\textrm{ess-liminf}_{y\to x}W(t,y)\le
W(t,x)$, then for every $\delta_j\to 0$,
$\delta_j>0$ there exists a sequence $(y_j^l)_j\in\N$ such that $y_j^l\to
x_l(t'')$, $W(t'',y_j^l)\le W(t'',x_l(t''))+\delta_j$ and $B^{y_j^l}$ is at
most countable.  Consider the following function defined on $[t',t'']$: 
\[
\fhi^l_j(t):=W(t,\zeta^{y_j^l}(t))+\int_{t'}^tL(s,\zeta^{y_j^l}(s),\omega )ds.
\] 
By the choice of $y_j^l$ and the hypotheses $\ref{i4})$, $\fhi_j^l$ is
differentiable a.e. with a nonnegative derivative. By the lower semicontinuity
of $W$ and the NDJ condition, it follows that $\fhi^l_j$ verifies the
hypotheses of Lemma \ref{lerv} and so $\fhi^l_j(t')\le\fhi^l_j(t'')$. Thus
\begin{equation}\label{7} 
W(t',\zeta^{y_j^l}(t'))\le
W(t'',\zeta^{y_j^l}(t''))+\int_{t'}^{t''}L(s,\zeta^{y_j^l}(s),\omega )ds.
\end{equation} 
Now, using the fact that $\zeta^{y_j^l}(t'')=y_j^l$ we obtain
\begin{eqnarray} 
W(t',\zeta^{y_j^l}(t')) & \le &
W(t'',y_j^l)+\int_{t'}^{t''}L(s,\zeta^{y_j^l}(s),\omega )ds\nonumber\\ & \le &
W(t'',x_l (t''))+\delta_j+\int_{t'}^{t''}L(s,\zeta^{y_j^l}(s),\omega
)ds.\label{7.1} 
\end{eqnarray} 
By Proposition 4.2, $\zeta^{y_j^l}(\cdot )\to
x_l(\cdot)$ as $j\to+\infty$ and so by the Lebesgue theorem and the lower
semicontinuity of $W$, passing to the limit as $j\to +\infty$ we obtain:
\begin{equation}\label{8} 
W(t',x_l (t'))\le W(t'',x_l
(t''))+\int_{t'}^{t''}L(s,x_l (s),\omega )ds. 
\end{equation} 
First consider the case $\left\{ (t,x_l (t)):t\in[t_0,T]\right\}\s Q$. Summing (\ref{8}) over each interval on which $u_l$ is constant we have
\begin{equation}\label{9}
W(t_0,x_l (t_0))\le W(T,x_l (T))+\int_{t_0}^TL(s,x_l (s),u_l (s))ds.
\end{equation}
Now, $x_l (T)=x^\ast (T)$ by definition and so, using (\ref{2}-\ref{4}) and \ref{i2})
\begin{eqnarray*}
W(t_0,x_l (t_0)) & \le & W(T,x^\ast (T))+\int_{t_0}^TL(s,x_l (s),u_l (s))ds\\
 & \le & \psi (T,x^\ast (T))+\int_{t_0}^TL(s,x_l (s),u_l (s))ds\nonumber\\
 & \le & V(t_0,x_0)+\frac{\eps}{2}-\int_{t_0}^TL(s,x^\ast (s),u^\ast (s))ds\nonumber\\
 & & +\int_{t_0}^TL(s,x_l (s),u_l (s))ds\nonumber\\
 & \le & V(t_0,x_0)+\eps < W(t_0,x_l (t_0))\nonumber.
\end{eqnarray*}
This is a contradiction.\\
Suppose now $\left\{ (t,x_l (t)):t\in[t_0,T]\right\}\not\s Q$. Define 
\begin{equation}\label{tau}
\hat\tau :=\inf\left\{ t\le T: (s,x_l (s))\in Q\quad\forall s\in [t,T]\right\}.
\end{equation}
In particular $(\hat \tau,x_l (\hat\tau ))\in\partial Q$. Using the same
argument to pass from (\ref{8}) to (\ref{9}), we obtain that for every
$\tau >\hat\tau$
\begin{equation}\label{bo}
W(\tau ,x_l (\tau))\le W(T,x^\ast (T))+\int_\tau^TL(s,x_l (s),u_l
(s))ds
\end{equation}
and so, using \ref{i2}) and (\ref{3})
\begin{eqnarray}
W(\tau ,x_l (\tau)) & \le & \psi(T,x^\ast (T))+\int_\tau^TL(s,x_l
(s),u_l (s))ds\nonumber\\
 & \le & V(t_0,x_0)+\frac{\eps}{2}-\int_{t_0}^TL(s,x^\ast (s),u^\ast
(s))ds\nonumber\\
 & & +\int_\tau^TL(s,x_l (s),u_l (s))ds.\label{lpos}
\end{eqnarray}
Using (\ref{1}), (\ref{4}) and \ref{i5}), we obtain for all $\tau >\hat\tau$
\begin{equation}
W(\tau ,x_l (\tau))\le V(t_0,x_0)+\eps\le W(t_0,x_0)-\eps.\label{Ref1}
\end{equation}
Passing to the liminf as $\tau\to\hat\tau$ and using the lower semicontinuity of $W$, we conclude
\begin{equation}\label{Ref2}
W(\hat\tau ,x_l (\hat\tau)) \le W(t_0,x_0)-\eps
\end{equation}
and so by $\ref{i3})$
\begin{equation}\label{Ref3}
W(t_0,x_0)\le \sup_{(t,x)\in Q} W(t,x)\le W(t_0,x_0)-\eps
\end{equation}
which is a contradiction. 

Now, we have to treat the case $V(t_0,x_0)=-\infty$. Since $W(t_0,x_0)>-\infty$ and $W$ is lower semicontinuous, we may find two constants $M>1$ and $\delta>0$ such that:
\[
W(t_0,x)>-M
\]
for every $x$ so that $\abs{x-x_0}<\delta$. Moreover we can find $u^\ast\in\U$ such that $x^\ast(\cdot):=x(\cdot;u^\ast,t_0,x_0)$ satisfies $(T,x^\ast(T))\in S$ and
\[
\int_{t_0}^TL(s,x^\ast(s),u^\ast(s))ds+\psi(T,x^\ast(T))\le -2M.
\]
With the same arguments of the first part of the proof we may find a control $u_l\in\U$ piecewise constant and left continuous such that, if $x_l(\cdot)$ is the trajectory $x(\cdot;u_l,T,x^\ast(T))$,
\[
\abs{x_l(t)-x^\ast(t)}<\frac{\delta}{2}\quad\forall t\in [t_0,T]
\]
and
\[
\abs{\int_{t_0}^T[L(s,x_l(s),u_l(s))-L(s,x^\ast(s),u^\ast(s))]ds}\le1.
\]
Repeating the same calculations as before, we obtain that
\begin{eqnarray*}
-M & \le & W(T,x^\ast(T))+\int_{t_0}^TL(s,x_l(s),u_l(s))ds\\
 & \le & \psi(T,x^\ast(T))+\int_{t_0}^TL(s,x^\ast(s),u^\ast(s))ds+1\\
 & \le & -2M+1
\end{eqnarray*}
which gives $M\le1$, a contradiction.

This concludes the proof of the theorem.
\end{proof}
\begin{corollary}
Let us suppose that $W$ satisfies all the hypotheses of the previous theorem. If moreover $W\ge V$ then $W=V$.
\end{corollary}

\begin{rmk}
If $W$ is produced by a synthesis procedure, the inequality $W\ge V$
always holds and so if $W$ satisfies all the hypotheses of Theorem
\ref{teo1} then $W$ coincides with the value function. \end{rmk}

Using the same techniques of the previous theorem, we can prove a corollary
for value functions generated by approximated syntheses, and give 
a bound of the error thus produced.

\begin{corollary}
Suppose (A1)-(A5) hold. Let $Q\s\Omega$ be an open subset containing $S$. Let $W:\overline Q\to\R$ be a lower semicontinuous function verifying the NDJ property along every time-varying vector field of the type $f(t,x,u)$ with $u\in U$ fixed. Moreover we assume that, for each $t$, ${\textrm{ess-liminf}}_{y\to x}W(t,y)\le W(t,x)$ and that there exist $\eps>0$ and $g\in L^1(\R)$, $g\ge0$, such that:
\renewcommand{\theenumi}{\roman{enumi}}
\renewcommand{\labelenumi}{\theenumi )}
\begin{enumerate}
\item $W\le\psi+\eps$ on $S$.
\item  At every point $(t,x)\in\partial Q$ one has
\[
W(t,x)=\sup_{(s,y)\in Q}W(s,y).
\]
\item There exists a countably $\mathcal H^n$-rectifiable set
 $A\s\Omega$ such that $W$ is differentiable on $Q\setminus A$ and satisfies
\[
W_s(s,y)+\inf_{\omega\in U}\left\{ W_y(s,y)\cdot
f(s,y,\omega)+L(s,y,\omega)\right\}\ge -\eps g(s)\quad \textrm{on }Q\setminus A.
\]
\item $L\ge -\eps g$.
\end{enumerate}
Then $W\le V+\eps(1+\| g\|_1)$ on $Q$.
\end{corollary}
\begin{proof}
Note that $L(t,x,u)+\eps g(t)\ge0$ and so
\begin{eqnarray*}
W(t_0,x_0) & \le & \!\!\!\!\!\!\!\!\!\!\!\!\inf_{\substack{u\in\U\vspace{3pt}\\ (T,x(T;u,t_0,x_0))\in S} }\!\!\!\left\{\int_{t_0}^T\!\!\!\!\!\!L(s,x(s;u,t_0,x_0),u(s))ds\!+\!\psi (T,x(T;u,t_0,x_0))\right\}\\
 & \le & V(t_0,x_0)+\eps(1+\| g\|_{L^1}).
\end{eqnarray*}
\end{proof}

\begin{rmk}
Notice that the value function of an optimal control problem has the NDJ
property along every possible direction as a consequence of the Dynamic
Programming Principle. Indeed, for every
$(t,y)\in\Omega\setminus S$ and for every admissible control $u\in\U$ (in
particular for every control $\omega\chi_I$, where $\omega\in U$ and $I$
bounded interval), the function 
\[
h\mapsto\int_t^{t+h}L(s,x(s;u,t,y),u(s))ds+V(t+h,x(t+h;u,t,y)) 
\]
is non
decreasing for $h\in[0,\delta]$ and $\delta$ small enough.

\noindent Instead, the hypothesis
\[
\textrm{ess-liminf}_{y\to x}W(t,y)\le W(t,x)
\]
for each $t$ fixed, says that, for every $\eps>0$ there exists a subset $V\s\{y\in\R^n:\abs{y-x}\le\eps\}$ of strictly positive Lebesgue measure such that
\[
\inf_{y\in V}W(t,y)\le W(t,x).
\]
So, if we consider a set $V_1\s\R^n$ of zero Lebesgue measure with $x$ as a cluster point, the set $V\setminus V_1$ has a strictly positive Lebesgue measure. 
In the proof of Theorem \ref{teo1} this fact is used to avoid the points $y$
for which $B^y$ is not countable. Moreover this hypothesis, coupled
with the lower semicontinuity of $W$, gives the following: \begin{itemize}
\item for each $t$, 
\[
W(t,x)=\liminf_{y\to x}W(t,y)=\textrm{ess-liminf}_{y\to x}W(t,y).
\]
\end{itemize}
\end{rmk}

%
%
\begin{rmk}\label{rip3}
Hypothesis \ref{i3}) of Theorem \ref{teo1} says that, in the case $Q\ne
\Omega$, the boundary of $Q$ must be a level set of the function $W$. We can
relax the same hypothesis in the following way:

\begin{itemize}
\item At every point $(t,x)\in\partial Q$ one has
\[
\liminf_{\substack{\tau\to t,y\to x\\(\tau,y)\in Q}}W(\tau,
y)\ge\sup_{(s,y)\in Q}W(s,y)
\]
\end{itemize}
and the conclusion of the theorem remains valid. Moreover if we define
with $R(t,x)$ the set of point reachable with an admissible control from
$(t,x)$, the previous condition can be replaced by
\[
\inf_{(s,y)\in R(t,x)\cap\partial Q}W(s,y)\ge W(t,x)
\]
and the conclusion still holds.
\end{rmk}
\\

The hypotheses of the positiveness of $L$ is almost optimal as the next
example shows. However, the Lagrangian $L$ may be negative on some 
region if trajectories 
can not 
stay for too long in such a region and 
one can relax the assumption \ref{i5}) as shown
in Remark \ref{RemL}4.
\vspace{0.5cm}\\
%
%
\begin{example}\label{ExaLpos}
\begin{figure}
\begin{center}
\includegraphics[height=6cm]{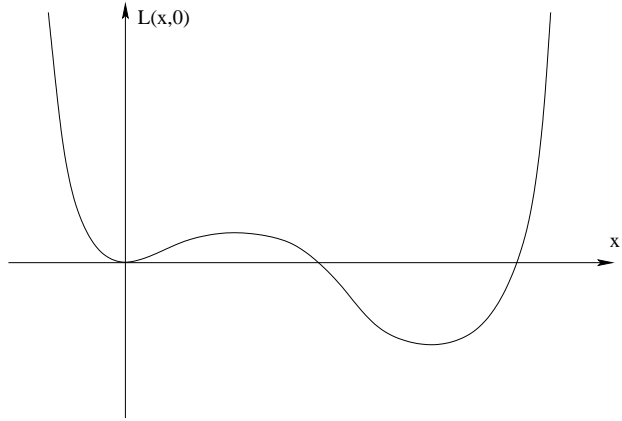}
\caption{$L(x,0)$ of Example 5.1.}
\label{fig}
\end{center}
\end{figure}
Consider the system $\dot x=u$, $U=[-1,1]\,$ and $\,\U=L^1(\R;U)$,
$\Omega=\R^2$, $S=\R\times\{ 0\}$, $Q=\R\times ]-1,1[$ with the Lagrangian
$L(t,x,u)=u^2+x^4-6x^3+7x^2$ (see Figure \ref{fig}) and $\psi\equiv 0$ on
$S$. Since the Lagrangian is negative in a region where the system can
stay for an arbitrary interval of times, clearly the value function for
this problem is equal to $-\infty$. If $W\equiv C$ on $\overline Q$ with
$C$ negative constant, then $W$ verifies all the hypotheses of the Theorem
\ref{teo1}, but \ref{i5}). In fact \ref{NU1}), \ref{i2}), \ref{i3})  are
obvious, while \ref{i4}) holds because $L$ is positive on $Q$ and $W$ is
differentiable on $Q$. So there exist infinitely many functions $W$ defined on
$\overline Q$ verifying the hypotheses of Theorem \ref{teo1}, but \ref{i5}),
which are not lower or equal to the value function $V$.
\end{example}\\

%
%
\renewcommand{\labelenumi}{\theenumi )}
\renewcommand{\theenumii}{\arabic{enumii}}
\renewcommand{\labelenumii}{\theenumi.\theenumii}
\begin{rmk}\label{RemL}
If one wants to eliminate hypothesis \ref{i5}) from the previous theorem, one
may assume one of the following conditions: \begin{enumerate}
\item Fix $\eps>0$ and $(\bar t,\bar x)\in Q$. We call $x_\eps:[\bar t,T]\to\R^n$ an $\eps$-quasi optimal trajectory ($\eps$-q.o.t.) for $(\bar t,\bar x)$ if:
\begin{enumerate}
\item $\exists\,u_\eps\in\U$ such that $\dot x_\eps(s)=f(s,x_\eps(s),u_\eps(s))$ for a.e. $s\in[\bar t,T]$,
\item $x_\eps(\bar t)=\bar x$,
\item $(T,x_\eps(T))\in S$,
\item $V(\bar t,\bar x)+\eps\ge\int_{\bar t}^TL(s,x_\eps(s),u_\eps(s))ds+\psi(T,x_\eps(T))$.
\end{enumerate}
Now define $Q_1$ as the set of point $(\bar t,\bar x)\in Q$ such that,
 for every $\eps>0$, there exists $x_\eps$, an $\eps$-q.o.t. for $(\bar
t,\bar x)$, satisfying $(s,x_\eps(s))\in Q$ for any $s\in[\bar t,T]$. What we
need is that $L\ge 0$ in $\Omega\setminus Q_1$. In fact, under this
assumption, we may suppose that $(s,x (s))\in\Omega\setminus Q_1$ for every
$s\in[t_0,\hat\tau [$, where $x$ is the trajectory defined in the proof of
Theorem \ref{teo1} and the time $\hat\tau$ is defined in (\ref{tau}). So the
integral $\int_{t_0}^{\hat\tau}L(s,x (s),u (s))ds$ is positive. Otherwise we
can assume $Q_1=Q$.
\item We can
also use an hypothesis similar to one given in \cite{ma}. For any $(\bar
t,\bar x)\in \Omega$ and $u\in\U$, let $x_{\bar t,\bar
x}(\cdot;u):=x(\cdot;u,\bar t,\bar x)$ be the solution to (\ref{cs})
associated to the control $u$. Consider the set $P$ consisting of those points
$(\bar t,\bar x)$ of $Q$ such that \[ \int_{\bar t}^TL(s,x_{\bar t,\bar
x}(s;u),u(s))ds\ge 0\qquad\forall T>t\quad\forall u\in\U. \] We have to
suppose that, if $(\bar t,\bar x)\in Q\setminus[P\cup S]$, there exist a
bounded and open set $B$, $(\bar t,\bar x)\in B\s Q$, $B\cap S=\emptyset$, so
that $\partial B\s Q$, and a positive number $M$ strictly less than \[
\inf_{u\in\U}\left\{ T>0:d((\bar t+T,x_{\bar t,\bar x}(\bar t+T;u)),\partial
B)\le d((\bar t,\bar x),\partial B)/2\right\} \] such that, for all $u\in\U$,
$(M+\bar t,x_{\bar t,\bar x}(M+\bar t;u))\in Q\cap P$ and  \[ \int_{\bar
t}^{\bar t+M}L(s,x_{\bar t,\bar x}(s),u(s))ds\ge 0, \] and this allow to
conclude the proof of Theorem \ref{teo1} without using $L\ge 0$ on the whole
space. 
\end{enumerate}\vspace{-.7cm} 
\end{rmk}

\begin{example}
Consider the system $\dot x=u$, $U=[-1,1]$, $\U=L^1(\R;U)$, $\Omega=\R^+\times\R$, $S=\R^+\times\{0\}$, $Q=\R^+\times ]-1,1[$, $\psi=0$ on $S$ and the Lagrangian defined by
\[
L(t,x,u):=
\left\{
\begin{array}{lrl}
u^2+x^2 & \textrm{if} & x\le 1\\
(u^2+1)(2-x)+(x-1)(u^2+Ct) & \textrm{if} & 1<x<2\\
u^2+x^2-6x+8+Ct & \textrm{if} & x\ge 2
\end{array}
\right.
\]
It is clear that this Lagrangian, for $C$ sufficiently big, 
satisfies the conditions a) and b) of the previous remark, even if it is not
positive outside $Q$. 
\end{example}

\begin{rmk}
We can relax hypotheses \ref{i3}) and \ref{i5}) with the following:
\begin{description}
\item[\textit{iii')}] the boundary $\partial Q$ is a level set of $W$;
\item[\textit{v')}] $L\ge0$ on $\Omega\setminus Q$.
\end{description} 
With these hypotheses, we can obtain an inequality of type (\ref{8}) for each
interval where the couple time-trajectory is in $Q$ and then, using $iii')$,
$v')$, the lower semicontinuity of $W$ and the NDJ property we can obtain
(\ref{9}).
\end{rmk}
%
%
%
%
\section{Problem with infinite time.}\label{se5}
\setcounter{equation}{0}
In this section we consider the control system (\ref{cs}) and assume that
(A-\ref{a1})-(A-5) hold with $0\leq C_R\leq C$ for some $C>0$
and every $R>0$.
Moreover we suppose that the target $S$ is a closed subset of
$\R\times\R^n$ which satisfies the structural property:
\begin{description}
\item[$(\ast)$] \hspace{.27cm}For any $T>0$, there exists $(t,x)\in S$ with $t\ge T$.
\end{description}
Let $S_1$ be an open neighborhood of $S$ contained in $\Omega$. Assume that the final cost $\psi$ is defined on $S_1$ and, if $d((t,x(t;u,t_0,x_0)),S)\to 0$ as $t\to+\infty$, then the trajectory $x(\cdot;u,t_0,x_0)$ is definitively in $S_1$, that is:
\begin{description}
\item[$(\ast\ast)$] \hspace{.3cm}$\!\!\!\exists\, T>t$ such that $(s,x(s;u,t_0,x_0))\in S_1$ for all $s\ge T$.
\end{description}
Define the value function:
\begin{equation}\label{vfit}
V(t_0,x_0):=\!\!\!\!\!\!\!\!\!\!\!\!\!\!\!\inf_{\substack{u\in\U\\[5pt]d((t,x(t;u,t_0,x_0)),S)\to 0\\ \text{as }t\to+\infty}}\!\!\!\!\!\!\!\!\!\!\!\!\!\left\{\!\int_{t_0}^{+\infty}\!\!\!\!\!\!\!\!\!\!\!L(s,x(s;u,t_0,x_0),u(s))ds\!+\!\limsup_{t\to+\infty}\psi(t,x(t;u,t_0,x_0))\right\}
\end{equation}
In other words, we consider only the trajectories that approach the target 
$S$ in infinite time. Notice that this condition does not imply that 
$(T,x(T))\not\in S$ for any $T\ge t_0$.\\

\begin{rmk}
The introduction of an open neighborhood of the target $S$ is due to a
technical reason and precisely to the fact that it is necessary to compare the
candidate value function to the final cost near the target. Notice that in the
following theorem the set $Q$ must contain $S_1$. For example we consider
$\Omega=\R^+\times\R$, $S=\R^+\times\{0\}$, $Q=\{(t,x):t>0,x<1/t\}$ and 
$S_1=\{(t,x):t>0,x<3/t\}$. If $(t,2/t)$, with $t>0$, is a trajectory,
then it is definitely in $S_1$, but it is never in $Q$.
\end{rmk}
\begin{theorem}\label{teo2}
Let $Q\s\Omega$ be an open subset containing $S_1$. 
Let $W:\overline Q\to\R$ be a lower semicontinuous function such that
\renewcommand{\theenumi}{\roman{enumi}}
\renewcommand{\labelenumi}{\theenumi )}
\begin{enumerate}
\item $W$ has the NDJ property along every time-varying vector field
of the type $f(t,x,u)$ with $u\in U$ fixed and for each $t$,
\[
\textrm{ess-liminf}_{y\to x}W(t,y)\le W(t,x).
\]
\item $W\le\psi$ on $S_1$.\label{ii2}
\item\label{ii3} At every point $(t,x)\in\partial Q$ one has
\[
W(t,x)=\sup_{(s,y)\in Q}W(s,y).
\]
\item\label{ii4} There exists a countable $\mathcal H^n$-rectifiable set $A\s\Omega$ such that $W$ is differentiable in $Q\setminus A$ and satisfies
\[
W_s(s,y)+\inf_{\omega\in U}\left\{ W_y(s,y)\cdot f(s,y,\omega)+L(s,y,\omega)\right\}\ge 0\quad \textrm{in }Q\setminus A.
\]
\item $L\ge 0$.\label{ii5} 
\end{enumerate}
Then $W\le V$ on $Q$. If $Q=\Omega$ we can drop hypotheses \ref{ii3}) and \ref{ii5}). 
\end{theorem}
\begin{proof}
Suppose by contradiction that there exists 
$(t_0,x_0)\in Q$ such that $W(t_0,x_0)>V(t_0,x_0)$. In particular $V(t_0,x_0)<+\infty$. First of all, let us consider the case $V(t_0,x_0)>-\infty$. As in the first part of the proof of Theorem \ref{teo1}, we can find $\eps>0$ and $\delta>0$ 
such that the following holds:
\begin{gather}
V(t_0,x_0)\le W(t_0,x_0)-2\eps\\
\abs{x-x_0}<\delta\quad\imp\quad W(t_0,x)>V(t_0,x_0)+\frac{3\eps}{2}\,\,.\label{21}
\end{gather}
We can choose $u^\ast\in\U$,
with the property that the trajectory $(t,x^\ast(t))$ approaches the target when $t\to+\infty$, and such that
\begin{equation}\label{eq2}
\int_{t_0}^{+\infty}L(s,x^\ast(s), u^\ast(s))ds+\limsup_{t\to+\infty}\psi(t,x^\ast (t))\le V(t_0,x_0)+\frac{\eps}{2}\,\,,
\end{equation}
where $x^\ast(\cdot)$ is the trajectory corresponding to the control $u^\ast$ such that $x^\ast(t_0)=x_0$.\\
Consider, now, a strictly increasing sequence of times $T_j>t_0$ 
converging to $+\infty$. We may suppose that $(t,x^\ast(t))\in Q$ for every $t\ge T_1$.
Fix $j\in\N$. For every $l\in\N$, there exists $u_j^l\in\U$ piecewise constant and left continuous such that $\| u_j^l-u^\ast\|_{L^p([t_0,T_j])}\le\frac{1}{l}$. So, by \cite[Th\'eor\`em IV.9]{brez}, we can extract a subsequence of $(u_j^l)_l$, denoted again with $(u_j^l)_l$, and we can find a function $h_j\in L^p([t_0,T_j])$ such that $\vert u_j^l\vert\le h_j$ a.e. for every $l\in\N$ and $u_j^l\to u^\ast$ for a.e. $t\in [t_0,T_j]$ as $l\to+\infty$.
Thus denoting
with $x_j^l(\cdot)$ the trajectory $x(\cdot;u_j^l,T_j,x^\ast(T_j))$, 
for $l$ sufficiently big we have (see Proposition 3.2)
\begin{equation}\label{lim1}
\abs{x_j^l(t)-x^\ast(t)}\le\frac{\delta}{2}\qquad\forall t\in [t_0,T_j]
\end{equation}
and then
\begin{equation}\label{eq1}
\abs{\int_{t_0}^{T_j}\left[ L(s,x_j^l(s),u_j^l(s))-L(s,x^\ast(s),u^\ast (s))\right]ds}\le\frac{\eps}{2}\,\,.
\end{equation}
Now, fix $l\in\N$ such that (\ref{lim1}) and (\ref{eq1}) hold. First, let us suppose that $\{ (t,x_j^l(t)):t\in[t_0,T_j]\}\s Q$. So, using Lemma \ref{fin}, Lemma \ref{lerv}, the same arguments as in the proof of Theorem \ref{teo1} and (\ref{eq1}) we conclude
\begin{eqnarray*}
W(t_0,x_j^l(t_0)) & \le & W(T_j,x_j^l(T_j))+\int_{t_0}^{T_j}L(s,x_j^l(s),u_j^l(s))ds\\
 & \le & W(T_j,x^\ast(T_j))+\int_{t_0}^{T_j}L(s,x^\ast(s),u^\ast(s))ds+\frac{\eps}{2}\,\,.
\end{eqnarray*}
Using (\ref{21}) and (\ref{lim1}) we have
\begin{equation}
V(t_0,x_0)+\frac{3\eps}{2}<W(T_j,x^\ast(T_j))+\int_{t_0}^{T_j}L(s,x^\ast(s),u^\ast(s))ds+\frac{\eps}{2}\,\,.
\end{equation}
Now consider the other case and precisely $\{ (t,x_j^l(t)):t\in[t_0,T_j]\}\not\s Q$. Define
\begin{equation}
\tau_j^l:=\inf\left\{ t\ge t_0:(s,x_j^l(s))\in Q\,\,\,\,\forall 
s\in[t,T_j]\right\}.
\end{equation}
Given $\tau_j^l<t<T_j$
\begin{equation}
W(t,x_j^l(t))\le W(T_j,x_j^l(T_j))+\int_t^{T_j}L(s,x_j^l(s),u_j^l(s))ds.
\end{equation}
Considering the fact that $(t,x_j^l(t))\to (\tau_j^l,x_j^l(\tau_j^l))$ 
as $t\to\tau_j^l$, $(\tau_j^l,x_j^l(\tau_j^l))\in\partial Q$ and (\ref{ii3}) we
obtain
\begin{equation}
W(t_0,x_0)\le 
W(T_j,x^\ast(T_j))+\int_{\tau_j^l}^{T_j}L(s,x_j^l(s),u_j^l(s))ds.
\end{equation}
We can now use the hypothesis $\ref{ii5})$, (\ref{21}) and (\ref{eq1}) in order to have
\begin{gather}
V(t_0,x_0)+\frac{3\eps}{2} <  W(t_0,x_0)\le\nonumber\\
 \le W(T_j,x^\ast(T_j))+\int_{t_0}^{T_j}L(s,x^\ast (s),u^\ast(s))ds+\frac{\eps}{2}\,\,.\label{29}
\end{gather}
In all cases we have that, for every $j\in\N$, 
\begin{equation}
V(t_0,x_0)+\frac{3\eps}{2} < W(T_j,x^\ast(T_j))+\int_{t_0}^{T_j}L(s,x^\ast (s),u^\ast(s))ds+\frac{\eps}{2}\,\,.
\end{equation}
So, applying the limsup as $j\to+\infty$ we get
\begin{eqnarray}
V(t_0,x_0)+\frac{3\eps}{2} & \le & \limsup_{j\to+\infty}W(T_j,x^\ast(T_j))+\int_{t_0}^{+\infty}\!\!\!L(s,x^\ast(s),u^\ast(s))ds+\frac{\eps}{2}\nonumber\\
& \le & \limsup_{t\to+\infty}W(t,x^\ast(t))+\int_{t_0}^{+\infty}\!\!\!L(s,x^\ast(s),u^\ast(s))ds+\frac{\eps}{2}\,\,.\nonumber
\end{eqnarray}
For $t$ sufficiently big, $(t,x^\ast(t))\in S_1$ and so, using (\ref{ii2}) and (\ref{eq2}),
\begin{eqnarray*}
V(t_0,x_0)+\frac{3\eps}{2} & \le & \limsup_{t\to+\infty}\psi(t,x^\ast(t))+\int_{t_0}^{+\infty}L(s,x^\ast(s),u^\ast(s))ds+\frac{\eps}{2}\\
 & \le & V(t_0,x_0)+\eps
\end{eqnarray*}
which implies
\begin{equation*}
V(t_0,x_0)\le V(t_0,x_0)-\frac{\eps}{2}
\end{equation*}
which is a contradiction.

It remains the case $V(t_0,x_0)=-\infty$. Since $W(t_0,x_0)>-\infty$ and $W$ is lower semicontinuous, we may find two constants $M>1$ and $\delta>0$ such that
\[
W(t_0,x)>-M
\]
for every $x$ so that $\abs{x-x_0}<\delta$. Moreover we can find $u^\ast\in\U$ such that $x^\ast(\cdot):=x(\cdot;u^\ast,t_0,x_0)$ approaches the target when $t\to+\infty$ and
\[
\int_{t_0}^{+\infty}L(s,x^\ast(s),u^\ast(s))ds+\limsup_{t\to+\infty}\psi(t,x^\ast(t))\le-2M.
\]
Consider a strictly increasing sequence of times $T_j>t_0$ converging to $+\infty$ and repeat the previous arguments in order to find a control $u_j^l\in\U$ piecewise constant, left continuous and such that, if $x_j^l(\cdot):=x(\cdot;u_j^l,T_j,x^\ast(T_j))$,
\[
\abs{x_j^l(t)-x^\ast(t)}\le\frac{\delta}{2}\quad\forall t\in[t_0,T_j]
\]
and
\[
\abs{\int_{t_0}^{T_j}[L(s,x_j^l(s),u_j^l(s))-L(s,x^\ast(s),u^\ast(s))]ds}\le1.
\]
Proceeding as before we obtain that
\begin{eqnarray*}
-M & \le & W(T_j,x^\ast(T_j))+\int_{t_0}^{T_j}L(s,x_j^l(s),u_j^l(s))ds\\
 & \le & W(T_j,x^\ast(T_j))+\int_{t_0}^{T_j}L(s,x^\ast(s),u^\ast(s))ds+1
\end{eqnarray*}
for every $j\in\N$. Passing to the limit we have:
\begin{eqnarray*}
-M & \le & \limsup_{j\to+\infty}W(T_j,x^\ast(T_j))+\int_{t_0}^{+\infty}L(s,x^\ast(s),u^\ast(s))ds+1\\
 & \le & \limsup_{t\to+\infty}W(t,x^\ast(t))+\int_{t_0}^{+\infty}L(s,x^\ast(s),u^\ast(s))ds+1\\
& \le & \limsup_{t\to+\infty}\psi(t,x^\ast(t))+\int_{t_0}^{+\infty}L(s,x^\ast(s),u^\ast(s))ds+1\\
& \le & -2M+1
\end{eqnarray*}
which gives $M\le 1$, a contradiction. 

So the theorem is proved.
\end{proof}
\begin{corollary}
Let $W$ satisfies all the hypotheses of the previous theorem and moreover $W\ge V$ where $V$ is defined in (\ref{vfit}). Then $W$ coincides with the value function.
\end{corollary}

\begin{rmk}
In theorem \ref{teo2} the condition \ref{ii2}) can be relaxed in the following way:
\[
\limsup_{t\to+\infty}W(t,x(t))\le\limsup_{t\to+\infty}\psi(t,x(t))
\]
for every $x(\cdot)$ solution to (\ref{cs}) such that $d((t,x(t)),S)\to0$ as $t\to+\infty$.

So, if one wants to minimize a Lagrangian cost without final cost, the condition becomes
\[
\limsup_{t\to+\infty}W(t,x(t))\le0
\]
for every $x(\cdot)$ with the above property.
\end{rmk}

%
%
\begin{rmk}
If we assume that there exists $\eta>0$ such that $S+B(0,\eta)\s S_1$, 
where $B(0,\eta)$ is the ball in $\R^{n+1}$ centered in $0$ with radius $\eta$,
then hypothesis $(\ast\ast)$ obviously holds. In fact suppose 
$d((t,x(t;u,t_0,x_0)),S)\to 0$ as $t\to+\infty$. Then there exists $T>0$ 
such 
that $d((s,x(s;u,t_0,x_0)),S)<\frac{\eta}{2}$ for all $s\ge T$. So we can choose an element $(t(s),y(s))\in S$ in order to have $d((s,x(s;u,t_0,x_0)),(t(s),y(s)))<\frac{\eta}{2}$ for all $s\ge T$. So the points $(s,x(s;u,t_0,x_0))\in S+B(0,\eta)\s S_1$ for every $s\ge T$.\vspace{.7cm}
\end{rmk}
%
%
\begin{rmk}
We obtain a generalization of Theorems \ref{teo1} and \ref{teo2} 
considering 
the same problem (\ref{cs}) with assumptions (A-\ref{a1})-(A-\ref{a4}), but we accept at the same time all the trajectories that hit the target in finite time or that tend to the target in infinite time. Obviously an analogous theorem as \ref{teo1} and \ref{teo2} holds.\vspace{.7cm}
\end{rmk}

%
%
\begin{rmk}
Also in this case we can substitute hypothesis $\ref{ii3})$ of Theorem 
\ref{teo2} in an analogous way as in Remark 5.3. Moreover we can eliminate
hypothesis \ref{ii5}) of Theorem \ref{teo2} in the same way as in Remark 
\ref{RemL}4.
\end{rmk}

%
%
%
%
\appendix
\renewcommand{\thesection}{\Alph{section}.}
\section{Viscosity solutions and value functions.}
\setcounter{equation}{0}
\setcounter{section}{1}
\setcounter{theorem}{0}
This appendix is intended to recall the notion of viscosity sub- and
super-solution and to state some known properties of the value function. 
Proof are analogous to those of \cite{bardi}.

Let $\Omega_1$ be an open subset of $\R\times\R^n$. We need the following definitions:
\begin{dfn}
Let $f:A\to\overline\R$ be a function where $A$ is an open subset of $\R^l$, for some $l\in\N\setminus\{ 0\}$. The lower semicontinuous envelope $f_\ast$ and the upper semicontinuous envelope $f^\ast$ of $f$ are defined by:
\begin{gather*}
f_\ast(x):=\lim_{r\to0^+}\inf\left\{ f(y):y\in A,\abs{y-x}\le r\right\},\\
f^\ast(x):=\lim_{r\to0^+}\sup\left\{ f(y):y\in A,\abs{y-x}\le r\right\}.
\end{gather*}
\end{dfn}
\begin{prop}
The lower semicontinuous (resp. upper semicontinuous) envelope of a function $f$ is a lower semicontinuous (resp. upper semicontinuous) function. More precisely, it is the greatest (resp. least) lower semicontinuous (resp. upper semicontinuous) function less or equal (resp. greater or equal) to $f$. Moreover $f$ is continuous if and only if $f_\ast=f^\ast$.
\end{prop}
\begin{dfn}
We say that a lower semicontinuous function $V:\Omega_1\to\overline\R$ is a viscosity super-solution to $F(t,x,D_tV,D_xV)=0$ in $\Omega_1$ if, for any $\fhi\in\mathcal C^1(\Omega_1)$ and for any $(t_0,x_0)\in\Omega_1$ point of local minimum for $V-\fhi$, one has $F^\ast(t_0,x_0,D_t\fhi(t_0,x_0),D_x\fhi(t_0,x_0))\ge 0$.
\end{dfn}
\begin{dfn}
We say that an upper semicontinuous function $V:\Omega_1\to\overline\R$ is a viscosity sub-solution to $F(t,x,D_tV,D_xV)=0$ in $\Omega_1$ if, for any $\fhi\in\mathcal C^1(\Omega_1)$ and for any $(t_0,x_0)\in\Omega_1$ point of local maximum for $V-\fhi$, one has $F_\ast (t_0,x_0,D_t\fhi(t_0,x_0),D_x\fhi(t_0,x_0))\le 0$.
\end{dfn} 

\begin{dfn}
We say that a function $V:\Omega_1\to\overline\R$ is a viscosity solution to $F(t,x,D_tV,D_xV)=0$ in $\Omega_1$ if $V_\ast$ is a viscosity super-solution and $V^\ast$ is a viscosity sub-solution to the equation.
\end{dfn}
\begin{rmk}
Note that the notion of viscosity solution is not bilateral, in the sense that the set of viscosity solution to $F=0$ and $-F=0$ in general are different.
\end{rmk}

Let us consider the following hypotheses:
\begin{description}
\item[(H-1)] The functions $f$ and $L$ are continuous in all the variables.
\item[(H-2)] $U$ is a bounded set.
\end{description}

We have the following:
\begin{prop}\label{pp1}
Let us assume (A-1)-(A-5) and (H-1)-(H-2).
Then the value function $V$ defined in
(\ref{vf}) satisfies the dynamic programming principle, that is
\[
V(t_0,x_0)=\!\!\!\!\!\!\!\inf_{\substack{u\in\U\vspace{3pt}\\(T,x(T;u,t_0,x_0))\in S}}\!\!\!\!\!\left\{\int_{t_0}^{T_1}\!\!\!\!\!\!L(s,x(s;u,t_0,x_0),u(s))ds\!+\!V(T_1,x(T_1;u,t_0,x_0))\right\}
\]
for every $(t_0,x_0)\in\Omega\setminus S$ and for every $T_1$ less than the minimum time to reach the target.
\end{prop}
An analogous proposition holds for the value function $V$ defined in (\ref{vfit}).\vspace{0.4cm}

%
%
Let us now state without proof the result that ensure that the value function is a viscosity solution to a Hamilton-Jacobi-Bellman equation.
\begin{theorem}\label{teoA}
Let us assume (A-1)-(A-5) and (H-1)-(H-2).
Then the value functions (\ref{vf}) and (\ref{vfit})
are viscosity solutions of 
\[
-V_s(t,x)-\inf_{\omega\in U}\left\{ f(t,x,\omega)\cdot V_y(t,x)+L(t,x,\omega)\right\}=0\quad \text{in }\Omega\setminus S.
\]
\end{theorem}

%
%
%
%
 
\end{document}